\documentclass[11pt]{article}

\usepackage[utf8]{inputenc}
\usepackage[english]{babel}
\usepackage[font=small]{caption}
\usepackage[T1]{fontenc}
\usepackage{amsmath,amssymb,amsfonts,amsthm}
\usepackage{mathrsfs} 
\usepackage{mathabx} 
\usepackage{bbm}
\usepackage{graphicx}
\usepackage{hyperref}
\usepackage[hmargin=2.5cm,vmargin=2cm]{geometry}
\usepackage[affil-it]{authblk}

\usepackage{enumitem}

\usepackage{color} 
\usepackage{longtable}

\theoremstyle{plain}
\newtheorem*{theorem*}{Theorem}
\newtheorem*{proposition*}{Proposition}
\newtheorem*{question*}{Question}

\theoremstyle{remark}

\def\N{\ensuremath{\mathbf{N}}}

\def\Z{\ensuremath{\mathbf{Z}}}

\def\ep{\varepsilon}

\def\E{\ensuremath{\mathbf{E}}}
\def\P{\ensuremath{\mathbf{P}}}

\def\Ind{\ensuremath{\mathbf{1}}}

\newcommand{\executeiffilenewer}[3]{%
\ifnum\pdfstrcmp{\pdffilemoddate{#1}}%
{\pdffilemoddate{#2}}>0%
{\immediate\write18{#3}}\fi%
}

\newcommand{\T}{\mathbb{T}}

\author{Itai Benjamini\thanks{e-mail: itai DOT benjamini AT weizmann DOT ac DOT il}, Pascal Maillard\thanks{e-mail: pascal DOT maillard AT weizmann DOT ac DOT il}}
\affil{The Weizmann Institute of Science\thanks{Department of Mathematics and Computer Science, Weizmann Institute of Science, POB~26, Rehovot~76100, Israel}}
\title{Point-to-point distance in first passage percolation on (tree)$\times \Z$}

\begin{document}

\maketitle

\begin{abstract}
\noindent We consider first passage percolation (FPP) on $\T_{d}\times G$, where $\T_d$ is the $d$-regular tree ($d\ge 3$) and $G$ is a graph containing an infinite ray $0,1,2,\ldots$. It is  shown that for a fixed vertex $v$ in the tree, the fluctuation of the distance in the FPP metric between the points $(v,0)$ and $(v,n)$ is of the order of at most $\log n$. We conjecture that the real fluctuations are of order $1$ and explain why.

\bigskip

\noindent \textbf{Keywords.} First passage percolation, Cayley graph, tree

\end{abstract}

\noindent Denote by $\T_d$ the $d$-regular tree (\(d\geq3)\), rooted at a vertex $\rho$. We consider FPP on $\T_{d}\times G$, where $G$ is a graph  containing an infinite ray $0,1,2,\ldots$ (for example, $G = \N$, $\Z$ or an infinite tree). That is, attach to each edge $e$ a random variable $X_e$, all $X_e$ being independent copies of a random variable $X\ge 0$ with $\E[X]<\infty$. For a path $\gamma$, denote by $|\gamma|$ the number of edges on the path and define $|\gamma|_X = \sum_{e\in \gamma} X_e$. Then define the random (pseudo-)metric
\[
 d_X(v,w) = \min\Big\{ |\gamma|_X\ \Big|\ \text{$\gamma$ is a path from $v$ to $w$}\Big\}.
\]
Write $D(n) = d_X((\rho,0),(\rho,n))$, i.e.\ the minimal distance between two points which are $n$ steps apart in the direction of the infinite ray.

We say hypothesis (H) is verified if
\begin{enumerate}
\item $\E[X^{1+\ep}]<\infty$ for some $\ep>0,$ and
\item there exist constants $C,K<\infty$, such that $\E[|\gamma_m|]<Cn^{K}$ for all $n$, where $\gamma_m$ is the path that minimizes $D(n)$ in  $\T_{d} \times G$.
\end{enumerate}
Note that 2. is verified for example if $X\ge c$ for some $c>0$, with $K=1$, because then $\E[|\gamma_m|] \le c^{-1}\E[|\gamma_m|_X] \le c^{-1}\E[X]n$.

\begin{theorem*}
Suppose hypothesis (H) is verified. Then, $(D(n) - \E[D(n)])/\log n$ is tight in $n$.
\end{theorem*}

We {\it conjecture} that  $D(n) - \E[D(n)]$ is tight (without rescaling). The rationale for this conjecture is that this is indeed the case for the graph $\T_{d-1,d}\times G$, where $\T_{d-1,d}$ is the rooted $d$-ary tree, i.e.\ the tree where the root has degree $d-1$ and every other vertex has degree $d$. This is the statement of the following proposition, which even does not need hypothesis (H):

\begin{proposition*}
\label{th:T}
 In $\T_{d-1,d}\times G$, the sequence $D(n) - \E[D(n)]$ is tight in $n$.
\end{proposition*}

First passage percolation is a model of random perturbation of a given geometry. It has mostly been studied in Euclidean space and lattices (see e.g.\ Howard \cite{Howard2004} for a review, although a bit outdated), and on trees, where it is also called the branching random walk \cite{Biggins2010,Shi2011,ZeitouniLectureNotes}. Other setups considered include the complete graph \cite{Janson1999,Bhamidi2012},   the Erd\"os-R\'enyi graph \cite{Bhamidi2011} and   a class of graphs admitting a certain recursive structure \cite{Benjamini2012} (see more on that below). However, to our knowledge, the present note is the first example where the fluctuation of the point-to-point distance in FPP on the Cayley graph of a finitely generated group (i.e.\ $\T_4\times\Z$) is shown to be small.

The graph $\T_d\times G$ can also be seen as an example of a ``large'' graph. The results of this note therefore add support to the common belief that the point-to-point distance of FPP in high-dimensional Euclidean space has small fluctuations.

We remark that the method  leading to the proof of the above theorem in general is not applicable for the study of the distance in the FPP metric between the points $(\rho,0)$ and $(v_n,0)$, for $v_n$ a vertex at distance $n$ from the root in $\T_d$. For example, in the case $G=\Z$, the minimizing path looks like a path in $\Z^2$ with additional ``handles'', and it is not clear for us whether the fluctuations are actually small (recall that in FPP on $\Z^2$, the fluctuations are believed to be of order $n^{1/3}$ (see e.g. \cite{Chatterjee2013a}) and up to now have only been proven to be of order at most $\sqrt{n/\log n}$ \cite{Benjamini2003}. We therefore ask the following question:

\begin{question*}
In $\T_d\times\Z$ or $\T_{d-1,d}\times\Z$, how big are the fluctuations of $d_X((\rho,0),(v_n,0))$, where $v_n$ is a vertex at distance $n$ from the root in the tree?
\end{question*}

We finally remark that even in $\T_d\times\T_d$, the current proof does not extend to the study to the FPP distance between two arbitrary vertices at distance $n$ apart in $\T_d\times\T_d$, although here we also conjecture that the fluctuations are of order 1.

\paragraph{Acknowledgement.} We thank an anonymous referee who has spotted some typographical errors. 

\section*{Proofs}
 The proof of the proposition uses a variant of an argument by Dekking and Host \cite{Dekking1991} on point-to-sphere distance in  FPP on a tree, which was generalized by Benjamini and Zeitouni \cite{Benjamini2012} to a large class of graphs, including $\T_{d-1,d}\times G$. For the point-to-sphere distance, the argument applies to every rooted graph $G$ containing two vertex-disjoint rooted subgraphs $G_1$ and $G_2$ which are isomorphic to $G$.\footnote{This is Property (1) in \cite{Benjamini2012}. Properties (2) and (3) are actually not needed, since on page 3 of that article, one can bound the right-hand side of the last inequality by $E\min(Z_n,Z_n')+KC$ and continue from that point on.} This argument can be adapted for the point-to-point distance in $\T_{d-1,d}\times G$ in the direction of the infinite ray in $G$. It fails for $\T_d\times G$, but not completely: It can be applied to an auxiliary graph, which ``almost'' looks like $\T_d\times G$. This however induces an error, which is the reason of the $\log n$ term appearing
in the statement of the theorem.

Before turning to the details, we introduce some more notation: Let $\T$ be a rooted tree and $v,w$ two vertices in $\T$. We say that $w$ is a \emph{descendant} of $v,$ if $v$ is contained in the direct path from the root to $w$.
We then denote by $\T|_v$ the subtree of $\T$ rooted at $v$, i.e.\ the subgraph of $\T$ spanned by the descendants of $v$, rooted at $v$.

\begin{proof}[Proof of the proposition]
Write for short $\T = \T_{d-1,d}$. Let $1$,$2$ denote two distinct children of the root in $\T$. For $i=1,2$, let $D_i(n)$ be the distance between $(i,0)$ and $(i,n)$ in the FPP metric restricted to the subgraph $\T|_i\times G$. Then $D_i(n)$ has the same law as $D(n)$. Furthermore, $D_1(n)$ and $D_2(n)$ are independent. Now, for $i=1,2$ and $j\in\N$, let $e_{i,j}$ be the edge between $(\rho,j)$ and $(i,j)$. Then
\[
 D(n) \le \min(D_1(n),D_2(n)) + X_{e_{1,0}} + X_{e_{1,n}} + X_{e_{2,0}} + X_{e_{2,n}}.
\]
Taking expectations and using the formula $\min(a,b) = (a+b)/2 - |a-b|/2$, we get
\[
 \E[D(n)] \le \tfrac 1 2 (\E[D_1(n)]+\E[D_2(n)]) + 4\E[X] - \tfrac 1 2 \E|D_1(n)-D_2(n)|.
\]
Since $\E[D_1(n)]=\E[D_2(n)]=\E[D(n)]$, this gives
\[
 \E|D_1(n)-D_2(n)| \le 8\E[X].
\]
Tightness follows from the inequality $\E|Z| \le \E|Z-Z'|$, which holds for any random variable $Z$ with $\E Z=0$ and with $Z'$ being an independent copy of $Z$.
\end{proof}

\begin{proof}[Proof of the theorem]
In the graph $\T_d\times G$ the situation is trickier: This graph does not contain two vertex-disjoint copies of itself. We will resolve this issue by considering an auxiliary graph first.

Write  $\T = \T_d$. Fix an integer $k$. Let $0$, $1$ and $2$ be three distinct children of the root in $\T$. For $i=0,1,2$, let $v_i$ be a vertex at distance $k$ from the root, which is a descendant of $i$. The auxiliary graph we consider is $\T'\times  G$, where $\T' = \T\backslash(\T|_{v_0})$, and the graph is rooted at $(\rho,0)$. In contrast to $\T \times G$, this graph does contain two vertex disjoint copies of itself, namely the graphs $\T|_1\times G$ and $\T|_2\times G$, rooted at $(v_1,0)$ and $(v_2,0)$, respectively.

Now let $D'=D'(n)$ be the FPP distance  between $(\rho,0)$ and $(\rho,n)$ in $\T'\times  G$ and $D_i'=D_i'(n)$ the FPP distance  between $(v_i,0)$ and $(v_i,n)$ in  $\T|_i\times G$, for $i=1,2$. Then note that $D'$, $D_1'$ and $D_2'$ have the same distribution and $D_1'$ and $D_2'$ are independent. Let $\gamma_{i,j}$ for $i=1,2$ and $j\in\N$ be the (unique) path from $(\rho,j)$ to $(v_i,j)$ in $\T\times\{j\}$ and let $|\gamma_{i,j}|_X$ be its length in the FPP metric. We then have
\[
 D' \le \min(D_1',D_2') + |\gamma_{1,0}|_X + |\gamma_{1,n}|_X + |\gamma_{2,0}|_X + |\gamma_{2,n}|_X,
\]
and taking expectations we get as in the previous proof
\[
  8\E[X]k \ge \E[|D_1'-D_2'|] \ge \E[|D'-\E[D']|].
\]
We now get back to the graph $\T\times G$. Let \(C\) denote some positive constant,
whose value may change from line to line. We claim that we can choose $k = k(n) = O(\log n)$, such that
\begin{equation}
\label{eq:DDprime}
 \E[|D-\E[D]|] \le \E[|D'-\E[D']|] + C.
\end{equation}
Together with the previous inequality, this will imply the statement of the theorem. To prove it, let $\gamma$ be the path from $(\rho,0)$ to $(\rho,n)$ of minimal length in the FPP metric on $\T\times G$ and let $V_\gamma$ be the projection on $\T$ of the set of vertices traversed by $\gamma$. Define the event $B = \{v_0\in V_\gamma\}$. Conditioned on $|V_\gamma|$, we have by symmetry,
\[
 \P\Big(B\,\Big|\,|V_\gamma|\Big) \le \frac{|V_\gamma|}{d\times(d-1)^{k-1}},
\]
since there are $d\times(d-1)^{k-1}$ vertices at  distance $k$ from the root in $\T_r$. Now, since $|V_\gamma|\le|\gamma|$, we have by hypothesis (H), with $k = \lceil\alpha\log_{d-1} n\rceil$, $\alpha > 0$,
\begin{equation}
 \label{eq:PB}
 \P(B) \le \E|V_\gamma|/(d\times(d-1)^{k-1}) \le Cn^{K}/(d-1)^k \le Cn^{K-\alpha}.
\end{equation}
Note that $D\le D'$ by definition, with $D = D'$ on the complement of  $B$. Together with the triangle inequality, this gives
\[
 \E|D-\E[D]|] \le \E[|D'-\E[D']|] + \E[|(D'-D)-\E[D'-D]|] \le C\log n + 2\E[D'\Ind_B].
\]
If $\gamma_0$ is the direct path from $(\rho,0)$ to $(\rho,n)$ along the ray $0,1,2,\ldots$, we have $D' \le |\gamma_0|_X$. Hypothesis (H) and Minkowski's inequality then give $\E[(D')^{1+\ep}] \le Cn^{1+\ep}$. Together with  H\"older's inequality and \eqref{eq:PB}, this yields the existence of $\alpha>0$, such that $\E[D'\Ind_B] < C$ for all $n$. This proves \eqref{eq:DDprime} and therefore finishes the proof of the theorem.
\end{proof}

\bibliography{fpp_tree_times_Z}

\end{document}